RESEARCH REPORTS

# REALIZATION SPACES OF 4-POLYTOPES ARE UNIVERSAL

JÜRGEN RICHTER-GEBERT AND GÜNTER M. ZIEGLER

Abstract. Let $P \subset \mathbb{R}^d$ be a $d$-dimensional polytope. The *realization space* of $P$ is the space of all polytopes $P' \subset \mathbb{R}^d$ that are combinatorially equivalent to $P$, modulo affine transformations. We report on work by the first author, which shows that realization spaces of 4-dimensional polytopes can be "arbitrarily bad": namely, for every primary semialgebraic set $V$ defined over $\mathbb{Z}$, there is a 4-polytope $P(V)$ whose realization space is "stably equivalent" to $V$. This implies that the realization space of a 4-polytope can have the homotopy type of an arbitrary finite simplicial complex, and that all algebraic numbers are needed to realize all 4-polytopes. The proof is constructive.

These results sharply contrast the 3-dimensional case, where realization spaces are contractible and all polytopes are realizable with integral coordinates (Steinitz's Theorem). No similar universality result was previously known in any fixed dimension.

## 1. Polytopes and their realization spaces

Polytopes have a long tradition as objects of mathematical study. Their historical roots reach back to the ancient Greek mathematicians, having a first highlight in the enumeration of the famous Platonic Solids. Already at this point strong impetus came from the fact that polytopes intimately relate topics from geometry, algebra and combinatorics. (The Platonic Solids solve a first enumerative question in polytopal geometry: to find all polytopes with a *flag transitive* symmetry group — a combinatorial concept.)

Received by the editors January 8, 1995.

1991 *Mathematics Subject Classification*. Primary 52B11, 52B40; Secondary 14P10, 51A25, 52B10, 52B30, 68Q15.

Supported by a DFG Gerhard-Hess-Forschungsförderungspreis.

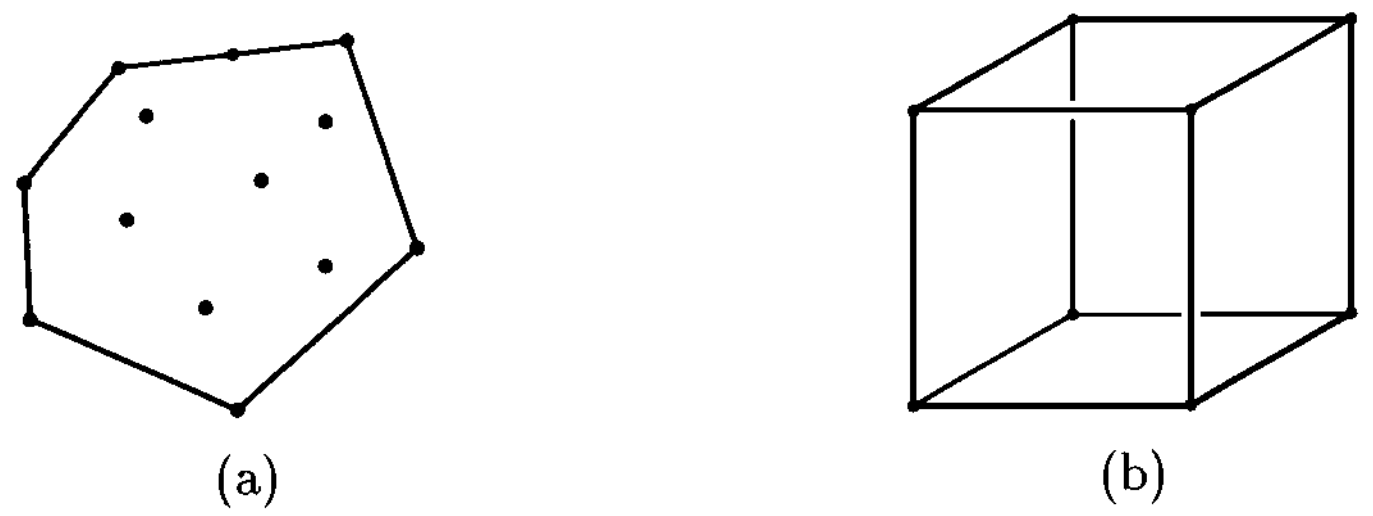

Figure 1

**Definition 1.** Let $\mathbf{P} = (\mathbf{p}_1, \ldots, \mathbf{p}_n) \in \mathbb{R}^{d\cdot n}$ be a finite collection of points that affinely span $\mathbb{R}^d$. The set

$$P = \mathbf{conv}(\mathbf{P}) := \Big\{ \sum_{i=1}^{n} \lambda_i \mathbf{p}_i \mid \sum_{i=1}^{n} \lambda_i = 1 \text{ and } \lambda_i \geq 0 \text{ for } i = 1, \ldots, n \Big\},$$

the *convex hull* of the point set $\mathbf{P}$, is called a *$d$-dimensional polytope* (a "*$d$-polytope*" for short). The *faces* of $P$ are the intersections $P \cap A$, where $A$ is an affine hyperplane that does not meet the relative interior of $P$. The *face lattice* of $P$ is the set of all faces of $P$, partially ordered by inclusion.

While a polytope is a geometric object, its face lattice is purely combinatorial in nature. Figure 1(a) illustrates a 2-polytope as the convex hull of finitely many points in the plane. We see that those points that are not in an extreme position have no contribution to the polytope itself. The points in extreme position (i.e., the 0-dimensional faces) are the *vertices* of a polytope. Figure 1(b) shows a cube as an example of a 3-dimensional polytope.

The need to structure the set of all polytopes of a fixed dimension leads to two main lines of study:

- to list all possible combinatorial types of polytopes (in other words, to decide which finite lattices correspond to polytopes and which do not),
- to describe the set of all realizations of a given combinatorial type.

The "set of all realizations" of a combinatorial type is formalized below by the concept of the *realization space* of a polytope. Besides their intrinsic importance for questions of real discrete geometry, such spaces appear in subjects as diverse as algebraic geometry (moduli spaces), differential topology (see Cairns's smoothing theory [8]), and nonlinear optimization (see Günzel et al. [10]).

Assume that in Definition 1 each point $\mathbf{p}_i$ for $i = 1, \ldots, n$ is a vertex of $P$. A *realization* of a polytope $P$ is a polytope $Q = \mathbf{conv}(\mathbf{q}_1, \ldots, \mathbf{q}_n)$ such that the face lattices of $P$ and $Q$ are isomorphic under the correspondence $\mathbf{p}_i \to \mathbf{q}_i$. The sequence of vertices $B = (\mathbf{p}_1, \ldots, \mathbf{p}_{d+1})$ is a *basis* of $\mathbf{P}$ if these points are affinely independent in any realization of $P$.

**Definition 2.** Let $P = \mathbf{conv}(\mathbf{p}_1, \ldots, \mathbf{p}_n) \subset \mathbb{R}^d$ be a $d$-polytope with $n$ vertices and with a basis $B = (\mathbf{p}_1, \ldots, \mathbf{p}_{d+1})$. The *realization space* $\mathcal{R}(P, B)$ is the set of all matrices $\mathbf{Q} = (\mathbf{q}_1, \ldots, \mathbf{q}_n) \in \mathbb{R}^{d\cdot n}$ for which $\mathbf{conv}(\mathbf{Q})$ is a realization of $P$ and $\mathbf{q}_i = \mathbf{p}_i$ for $i = 1, \ldots, d+1$.

It turns out that the realization space $\mathcal{R}(P, B)$ is essentially (up to "stable equivalence," see below) independent of the choice of a basis. Hence it makes sense to speak of *the* realization space $\mathcal{R}(P)$ of a polytope.

Every realization space is a *primary semialgebraic set* defined over $\mathbb{Z}$: it is the set of solutions of a finite system of polynomial equations $f_i(x) = 0$ and strict inequalities $g_j(x) > 0$, where the $f_i$ and $g_j$ are polynomials with integer coefficients on some $\mathbb{R}^N$. To see this, one checks that the realization space is the set of all matrices $\mathbf{Q} \in \mathbb{R}^{d \cdot n}$ for which some entries are fixed, and the determinants of certain $d \times d$ minors have to be positive, or negative, or zero. For a general *semialgebraic set* one also admits non-strict inequalities $h_k \geq 0$. Thus, for example, the set $\{0, 1\}$ and the open interval $]0, 1[\subset \mathbb{R}$ are primary semialgebraic sets, while the closed interval $[0, 1]$ is a semialgebraic set in $\mathbb{R}$ that is not primary.

In this research report we present a Universality Theorem proved by the first author [17], stating that all primary semialgebraic sets are in a suitable sense "stably equivalent" (see Section 3) to the realization spaces of suitable 4-polytopes.

## 2. Old and new results on realization spaces

What does the realization space of a polytope look like? Which algebraic numbers are needed to coordinatize the vertex set of a given $d$-dimensional polytope? How can one tell whether a finite lattice is the face lattice of a polytope or not?

For 3-dimensional polytopes, Steinitz's work [19, 20] answered these basic questions about realization spaces more then seventy years ago. In particular, Steinitz's "Fundamentalsatz der konvexen Typen" (today known as *Steinitz's Theorem*) and its modern relatives (see [9] and [23]) provide complete answers to these questions for this special case.

**Steinitz's Theorem, 1922.** *A graph $G$ is the edge graph of a 3-polytope if and only if $G$ is simple, planar and 3-connected.*

Here *simple* means that $G$ contains no "parallel edges" and no "loops". A graph $G$ is *planar* if it can be drawn in the plane without intersections, and it is 3-*connected* if between any two vertices there are three disjoint paths in $G$.

As corollaries and by inspection of the proof of Steinitz's Theorem one obtains:

- For every 3-polytope $P \subseteq \mathbb{R}^3$ the realization space $\mathcal{R}(P)$ is a smooth open ball. (This ball has dimension $e - 6$, if $P$ has $e$ edges.)
- For every 3-polytope $P$ the space $\mathcal{R}(P)$ contains rational points, that is, every 3-polytope can be realized with integral vertex coordinates.
- The shape of one 2-face in the boundary of a 3-polytope $P$ can be arbitrarily prescribed, that is, the canonical map $\mathcal{R}(P) \to \mathcal{R}(F)$ is surjective for every facet $F \subseteq P$ (Barnette & Grünbaum [2]).

Similar statements for $d$-polytopes that have at most $d + 3$ vertices were proved by Mani [14] and Kleinschmidt [13].

Over the years, it became increasingly clear that no similar positive answer could be expected for high-dimensional polytopes. Duality theory ("Gale diagrams" [9, 23]) was used to construct a non-rational 8-polytope with 12 vertices (Perles, 1967 [9]). Also Mnëv's famous Universality Theorem for oriented matroids [15, 16, 5, 10] via Gale diagrams implies a universality theorem for $d$-polytopes with $d + 4$ vertices: in general for such polytopes the realization spaces can be arbitrarily complicated.

**Mnëv's Universality Theorem for Polytopes, 1986.** *For every primary semialgebraic set $V$ defined over $\mathbb{Z}$ there is some $d \in \mathbb{N}$ and a $d$-polytope $P$ with $d+4$ vertices whose realization space is stably equivalent to $V$.*

In particular, stable equivalence implies homotopy equivalence. Stable equivalence also preserves "algebraic complexity" (see Lemma 4 below). As a consequence, all algebraic numbers are needed to coordinatize all $d$-polytopes with $d+4$ vertices.

However, Gale diagram techniques and their variants (such as "Lawrence extensions", see Section 4) do not provide systematic construction methods for $d$-polytopes in any fixed dimension $d$. Only a "sporadic" example of a 4-polytope with disconnected realization space was constructed (Bokowski, Ewald, and Kleinschmidt [6, 15, 7]). More sporadic examples showed that the shapes of 3-faces of 4-polytopes (Kleinschmidt [12], Barnette [1]) and of 2-faces of 5-polytopes (Ziegler [23]) *cannot* be prescribed arbitrarily. Until now no general construction techniques to produce polytopes with controllably bad behavior for any fixed dimension $d$ were known. The $\sigma$-construction presented in [21] for that purpose turned out to be incorrect [23].

In the following we report on the first author's recent work [17] that produces a complete systematic Universality Theorem for polytopes of dimension 4 (and thus for $d$-polytopes of any fixed dimension $d \geq 4$). In this work, all the "sporadic" examples of 4- and 5-polytopes are explained in terms of Lawrence extensions of planar and 3-dimensional point configurations (see Section 4). Sporadic examples (such as a new 4-polytope for which the shape of a hexagonal 2-face cannot be prescribed) are used as the "basic building blocks" of the construction for the following result [17].

**Main Theorem** (Richter-Gebert, 1994). *For every primary semi-algebraic set $V$ defined over $\mathbb{Z}$ there is a 4-polytope $P(V)$ whose realization space $\mathcal{R}(P(V))$ is stably equivalent to $V$. Moreover, the face lattice of $P(V)$ can be generated from the defining equations of $V$ in polynomial time.*

**Corollaries of the Main Theorem and its construction techniques.**

- *There is a non-rational 4-polytope with 34 vertices.*
- *All algebraic numbers are needed to coordinatize all 4-polytopes.*
- *The realizability problem for 4-polytopes is NP-hard.*
- *The realizability problem for 4-polytopes is (polynomial time) equivalent to the "Existential Theory of the Reals".*
- *For every finite simplicial complex $\Delta$, there is a 4-polytope whose realization space is homotopy equivalent to $\Delta$.*
- *There is a 4-polytope for which the shape of a 2-face cannot be arbitrarily prescribed.*
- *Boundary complexes of 4-polytopes cannot be characterized by excluding a finite set of "forbidden minors".*
- *In order to realize all combinatorial types of rational $n$-vertex 4-polytopes with vertex coordinates in $\{0, 1, \ldots, f(n)\}$, the "coordinate size" function $f(n)$ has to be at least a doubly exponential function in $n$.*

In particular these implications solve all the problems that were recently emphasized in "*Three problems about* 4*-polytopes*" [22]. They also solve [23, Problems 5.11*, 6.10*, and 6.11*].

## 3. Stable equivalence

The concept we use to compare realization spaces with general primary semialgebraic sets is *stable equivalence*. Although such a concept has been used by different authors, the precise definitions they used (see [11, 10, 15, 16, 18]) vary substantially in their technical content. The common idea is that semialgebraic sets that only differ by a "trivial fibration" and a rational change of coordinates should be considered as stably equivalent, while semialgebraic sets that differ in certain "characteristic properties" should turn out not to be stably equivalent. In particular, stable equivalence should preserve the homotopy type and respect the algebraic complexity and the singularity structure. We now present a concept of stable equivalence that is stronger than all previously used notions.

Let $V \subseteq \mathbb{R}^n$ and $W \subseteq \mathbb{R}^{n+d}$ be semialgebraic sets with $\pi(W) = V$, where $\pi$ is the canonical projection $\pi : \mathbb{R}^{n+d} \to \mathbb{R}^n$ that deletes the last $d$ coordinates. $V$ is a *stable projection* of $W$ if $W$ has the form

$$W = \Big\{(\mathbf{v}\,,\,\mathbf{v}') \in \mathbb{R}^{n+d} \mid \mathbf{v} \in V \text{ and } \phi_i(\mathbf{v}) \cdot \mathbf{v}' > 0\,;\ \psi_j(\mathbf{v}) \cdot \mathbf{v}' = 0 \quad \text{for all } i \in X\,;\ j \in Y\Big\}.$$

Here $X$ and $Y$ denote finite (possibly empty) index sets. For $i \in X$ and $j \in Y$ the functions $\phi_i$ and $\psi_j$ have to be polynomial functions

$$\phi_i = (\phi_i^1\,, \ldots\,, \phi_i^d) : \mathbb{R}^n \to (\mathbb{R}^d)^* \text{ with } \phi_i^k \in \mathbb{Z}[x_1\,, \ldots\,, x_n]$$

and

$$\psi_j = (\psi_j^1\,, \ldots\,, \psi_j^d) : \mathbb{R}^n \to (\mathbb{R}^d)^* \text{ with } \psi_j^k \in \mathbb{Z}[x_1\,, \ldots\,, x_n]\,,$$

that associate to every element of $\mathbb{R}^n$ a linear functional on $\mathbb{R}^d$.

Two semialgebraic sets $V$ and $W$ are *rationally equivalent* if there exists a homeomorphism $f : V \to W$ such that both $f$ and $f^{-1}$ are rational functions.

**Definition 3.** Two semialgebraic sets $V$ and $W$ are *stably equivalent*, denoted $V \approx W$, if they are in the same equivalence class with respect to the equivalence relation generated by stable projections and rational equivalence.

**Lemma 1.** *Let $V \approx W$, with $V \subseteq \mathbb{R}^n$ and $W \subseteq \mathbb{R}^m$, be a pair of stably equivalent semialgebraic sets, and let $A$ be a subfield of the algebraic numbers. We have*

(i) *$V$ and $W$ are homotopy equivalent.*
(ii) $V \cap A^n = \emptyset \iff W \cap A^m = \emptyset$.
(iii) *$V$ and $W$ have equivalent "singularity structure".*

## 4. Polytopal tools

Lawrence extensions and connected sums are elementary geometric operations on polytopes that form the basis for the constructions presented in [17]. They are very simple and innocent-looking operations, but still very powerful.

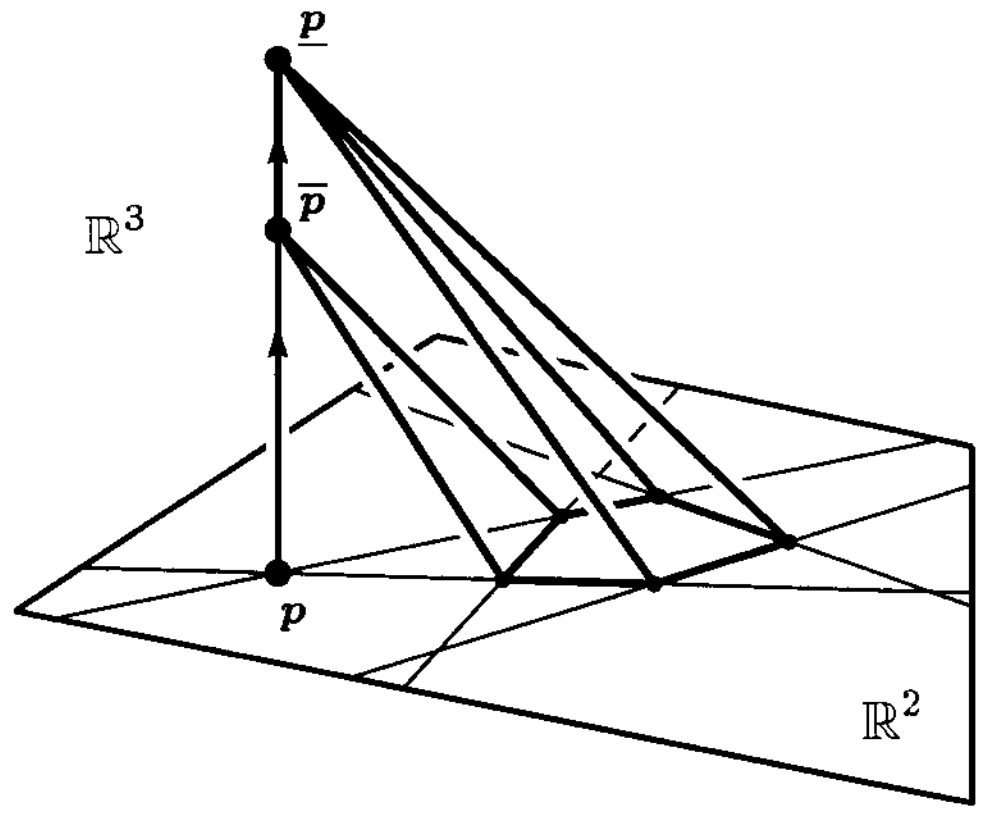

Figure 2

For *Lawrence extensions* the basic operation is the following: take a point $\mathbf{p}$ in a $d$-dimensional point configuration, and replace it by two new points $\overline{\mathbf{p}}$ and $\underline{\mathbf{p}}$ that lie on a ray that starts at the original point but goes off in some "new" direction into $(d+1)$-dimensional space (see Figure 2).

Every such Lawrence extension increases both the dimension of a point configuration, and its number of points, by $1$. Note that although the original point is deleted in the construction, it is still implicitly present: it can be "reconstructed" as the intersection of the line spanned by the two new points with the $d$-hyperplane spanned by the original point configuration.

The "classical" use of Lawrence extensions [4, 3, 15] starts with a 2-dimensional configuration of $n$ points and performs Lawrence extensions on all these points, one after the other.

The resulting configuration of $2n$ points is the vertex set of an $(n+2)$-dimensional polytope, the *Lawrence polytope* of the point configuration.

Every realization of the Lawrence polytope determines a realization of the original point configuration, including all collinearities and all orientations of triples. In fact, the realization spaces of the Lawrence polytope and the planar configuration can be shown to be stably equivalent. This can be used to lift Mnëv's Universality Theorem from planar point configurations (oriented matroids) to $d$-polytopes.

If one wants to stay within the realm of 4-polytopes, then it is not permissible to use more than two Lawrence extensions. However, careful use of just one or two Lawrence extensions on some points outside a 2- or 3-polytope leads to extremely interesting and useful polytopes — among them are the basic building blocks for the Main Theorem. Here space permits us only to sketch two such examples from [17].

First consider Pascal's theorem: if the vertices of a hexagon $H$ lie on an ellipse, then the intersection points of the lines spanned by opposite edges are collinear. If one takes the vertices of the hexagon together with the three intersection points as the initial configuration and performs Lawrence extensions on the three intersection points, then one obtains a 5-dimensional polytope with 12 vertices that has the original hexagon as a 2-face (Figure 3(a)). Furthermore, the six new points that are generated by the Lawrence extensions are affinely dependent and form a facet $F$ of the resulting 5-polytope.

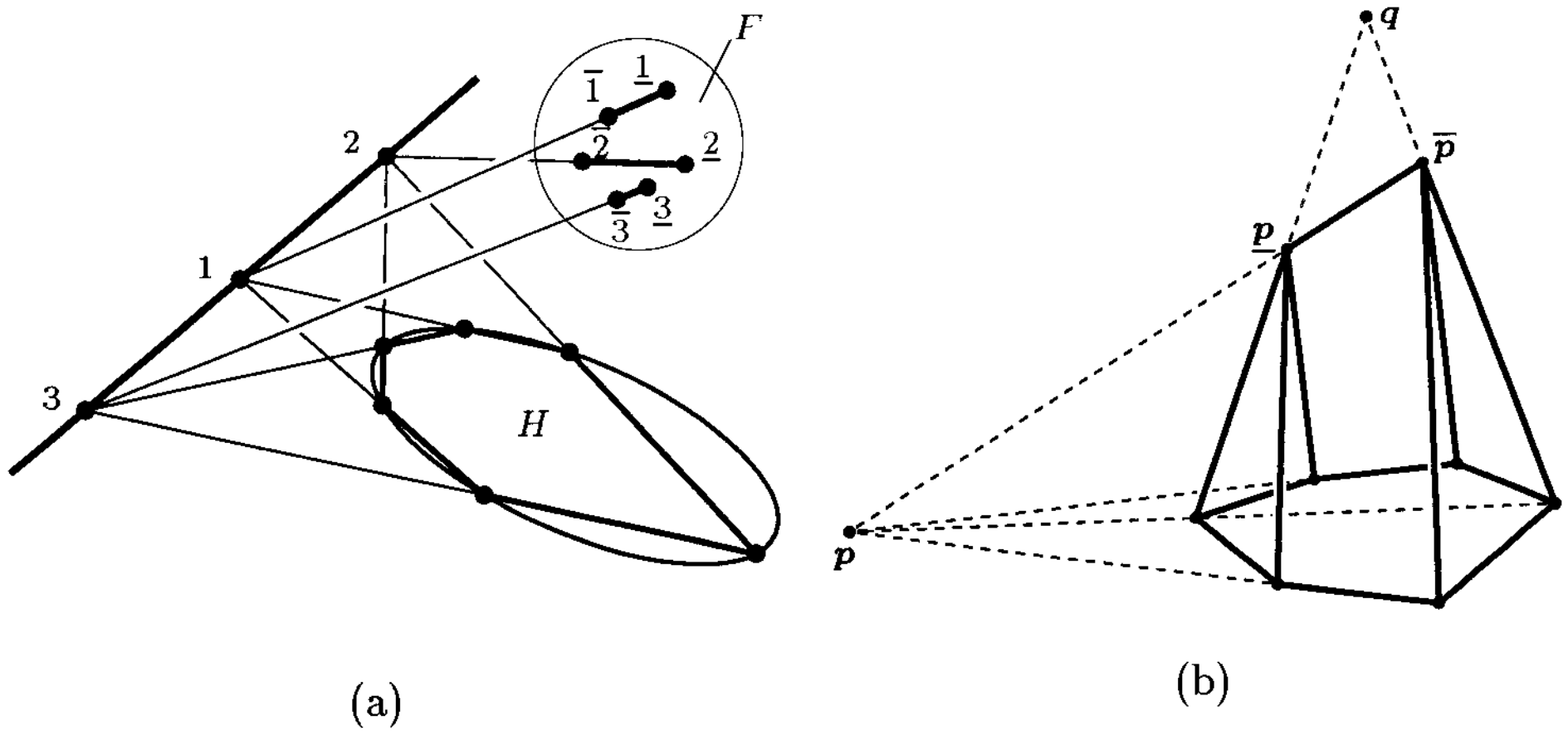

Figure 3

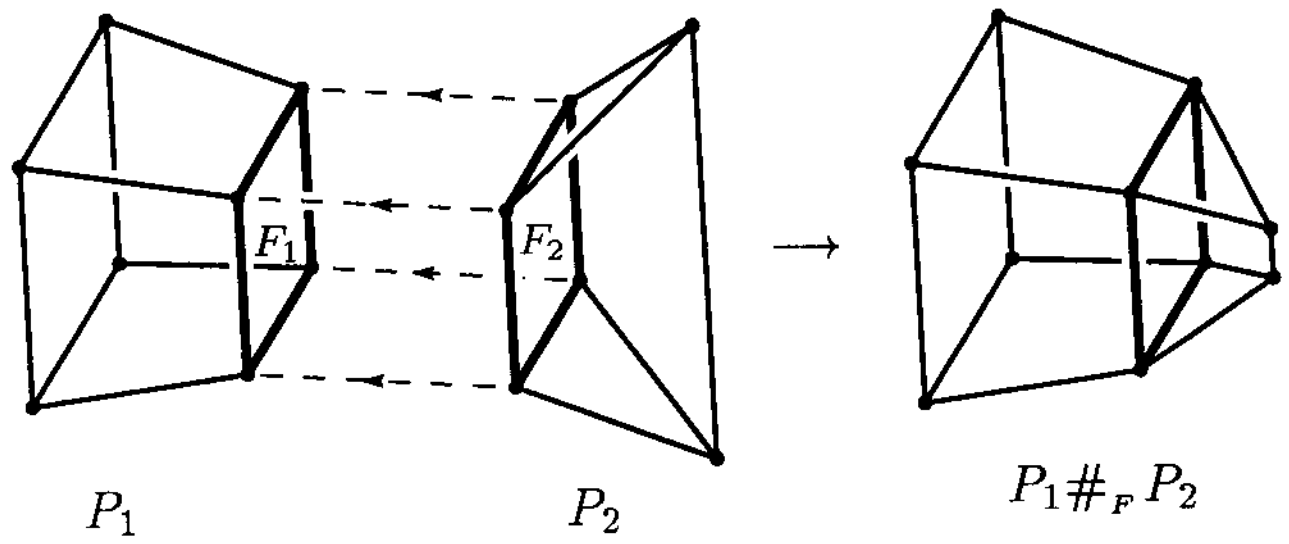

Figure 4

The combinatorics of Lawrence extensions, together with Pascal's Theorem, imply that in *every* realization of this 5-polytope, the hexagon has to have its vertices on an ellipse. (This is the 5-polytope constructed in [23, p. 175].)

This is quite surprising — see the corresponding corollary of Steinitz's Theorem for 3-polytopes. However, what was needed for the Main Theorem was a 4-dimensional polytope with the same property. The search for such an example was an open problem in [23, Problem 5.11]. It is now solved by performing a Lawrence extension on the point $\mathbf{q}$ in Figure 3(b), outside a "tent" over a hexagon (a 3-polytope that arises by performing a Lawrence extension on a certain point outside a hexagon). We leave the proof for this example to the reader (see [17]) and proceed from Lawrence extensions to the other basic construction method.

*Connected sums* of polytopes are obtained as follows. Assume that one is given two $d$-polytopes, $P_1$ and $P_2$, that have projectively equivalent facets $F_1$ and $F_2$. We use $F$ to denote the combinatorial type of $F_1 \cong F_2$. Then, using a projective transformation, one can "merge" $P_1$ and $P_2$ into a more complicated polytope, the *connected sum* $Q := P_1 \#_F P_2$. The polytope $Q$ has all the facets of $P_1$ and $P_2$, except for $F_1$ and $F_2$. However, the boundary complex $\partial F$, consisting of all the proper faces of $F$, is still present in $Q$ (Figure 4).

Now, if one takes an arbitrary realization of $Q$, then it is not in general true that this realization arises as a connected sum of realizations of $P_1$ and of $P_2$: in a "bad" realization of $Q$ the boundary complex $\partial F$ may not be flat. In fact, in dimension $d = 3$ one can see that the complex $\partial F$ in $Q$ is necessarily flat if and only if $F$ is a triangle facet. In dimension 4, there are several other types of facets that are "necessarily flat", among them are pyramids, prisms, and tents. Only such necessarily flat facets are used in connected sum operations for the proof of the Main Theorem.

## 5. Sketch of the proof

The proof of the Main Theorem is constructive. It starts with the defining equations of a semialgebraic set and uses them to explicitly construct the face lattice of a 4-polytope. A result of Shor [18] is used, which states that every primary semialgebraic set $V$ is stably equivalent to a semialgebraic set $V' \in \mathbb{R}^n$ whose defining inequalities fix a total order

$$1 = x_1 < x_2 < x_3 < \ldots < x_n$$

on the variables and for which all defining equations have the form

$$x_i + x_j = x_k \quad \text{or} \quad x_i \cdot x_j = x_k$$

for certain $1 \le i \le j < k \le n$. Such a set of defining equations and inequalities is a *Shor normal form* of $V$. Thus only elementary addition and multiplication have to be modelled: they are encoded into certain polytopes for which certain 2-faces are not prescribable.

In the following we briefly describe how for a given primary semialgebraic set $V$ in Shor normal form a corresponding polytope $P(V)$ can be constructed whose realization space is stably equivalent to $V$. While Lawrence extensions are used to generate "basic building blocks", the connected sum operation is used to combine these blocks to larger semantic units.

(i) The initial building blocks generated by Lawrence extensions are
- a 4-polytope that contains a hexagonal 2-face with vertices $1, \ldots, 6$, in this order, such that in every realization of $X$ the lines $1 \vee 4$, $2 \vee 3$ and $5 \vee 6$, are concurrent (see Section 4, Figure 3(b));
- "forgetful transmitter" polytopes $T_n$ that contain as faces an $n$-gon $G$ and an $(n-1)$-gon $G'$, such that in every realization of $T_n$ the configuration of the lines that are supported by the edges of $G'$ is projectively equivalent to the configuration of lines that is determined by certain $n-1$ edges of $G$;
- "connector" polytopes $C_n$ that contain three $n$-gons $G_1$, $G_2$ and $G_3$ that are projectively equivalent to each other in every realization of $C_n$.

(ii) Taking connected sums of the basic building blocks in (i), we obtain polytopes $P^+$ and $P^\times$ that model addition and multiplication. These two polytopes both contain 12-gons $G$ with edges labelled by

$$0, 1, i, j, k, \infty, 0', 1', i', j', k', \infty'$$

in this order. In each realization of $P^+$ or $P^\times$ the six intersections $\alpha^* = \alpha \cap \alpha'$ of opposite edge supporting lines of $G$ lie on a line $\ell$. The points $0^*$, $1^*$ and $\infty^*$ define a projective scale $\sigma$ on $\ell$. Furthermore, $P^+$ (resp.

$P^{\times}$) is realizable if and only if $\sigma(i^*)+\sigma(j^*)=\sigma(k^*)$ (resp. $\sigma(i^*)\cdot\sigma(j^*)=\sigma(k^*)$). (Special care has to be taken in the case $i=j$.)

(iii) Again by connected sum operations these addition and multiplication polytopes are composed to obtain the polytope $P(V)$ that contains a $(2n+6)$-gon $G$ with edges labelled by

$$0\,,\,1\,,\,2\,,\,\dots\,,\,n\,,\,\infty\,,\,0'\,,\,1'\,,\,2'\,,\,\dots\,,\,n'\,,\,\infty'$$

in this order. In each realization of $P(V)$ the $n+3$ intersections $\alpha^*=\alpha\cap\alpha'$ of opposite edge supporting lines of $G$ lie on a line $\ell$. The points $0^*$, $1^*$ and $\infty^*$ define a projective scale $\sigma$ on $\ell$. Addition and multiplication polytopes are adjoined according to the defining equations of $V$. In this way the points of $V$ are in one-to-one correspondence with the values $\sigma(1^*),\dots,\sigma(n^*)$ in all possible realizations of $P(V)$.

Thus, modulo projective equivalence, $P(V)$ contains a centrally symmetric $(2n+6)$-gon whose slopes of opposite edges in any realization of $P(V)$ encode the coordinates of the corresponding point in the semialgebraic set $V$:

$$\mathcal{R}(P(V)) \;\approx\; \{(\sigma(1^*),\dots,\sigma(n^*)) \mid P\in\mathcal{R}(P(V))\} \;=\; V' \;\approx\; V.$$

## References

[1] D.W. Barnette, *Two "simple" 3-spheres,* Discrete Math. **67** (1987), 97–99.

[2] D. Barnette and B. Grünbaum, *Preassigning the shape of a face,* Pacific J. Math. **32** (1970), 299–302.

[3] M. Bayer and B. Sturmfels, *Lawrence polytopes,* Canad. J. Math. **42** (1990), 62–79.

[4] L.J. Billera and B.S. Munson, *Polarity and inner products in oriented matroids,* European J. Combin. **5** (1984), 293–308.

[5] A. Björner, M. Las Vergnas, B. Sturmfels, N. White, and G.M. Ziegler, *Oriented matroids,* Encyclopedia of Mathematics, vol. 46, Cambridge Univ. Press, London and New York, 1993.

[6] J. Bokowski, G. Ewald, and P. Kleinschmidt, *On combinatorial and affine automorphisms of polytopes,* Israel J. Math. **47** (1984), 123–130.

[7] J. Bokowski and A. Guedes de Oliveira, *Simplicial convex 4-polytopes do not have the isotopy property,* Portugal. Math. **47** (1990), 309–318.

[8] S.S. Cairns, *Homeomorphisms between topological manifolds and analytic manifolds,* Ann. of Math. **41** (1940), 796–808.

[9] B. Grünbaum, *Convex polytopes,* Interscience, London, 1967; revised edition (V. Klee and P. Kleinschmidt, eds.), Springer-Verlag (in preparation).

[10] H. Günzel, *The universal partition theorem for oriented matroids,* Discrete Comput. Geom. (to appear).

[11] H. Günzel, R. Hirabayashi, and H. Th. Jongen, *Multiparametric optimization: On stable singularities occurring in combinatorial partition codes,* Control Cybernet. **22** (1994), 153–167.

[12] P. Kleinschmidt, *On facets with non-arbitrary shapes,* Pacific J. Math. **65** (1976), 511–515.

[13] ———, *Sphären mit wenigen Ecken,* Geom. Dedicata **5** (1976), 97–101.

[14] P. Mani, *Spheres with few vertices,* J. Combin. Theory Ser. A **13** (1972), 346–352.

[15] N.E. Mnëv, *The universality theorems on the classification problem of configuration varieties and convex polytopes varieties,* Topology and Geometry – Rohlin Seminar (O. Ya. Viro, ed.), Lecture Notes in Math., vol. 1346, Springer-Verlag, Heidelberg, 1988, pp. 527–544.

[16] ———, *The universality theorems on the oriented matroid stratification of the space of real matrices,* Applied Geometry and Discrete Mathematics – The Victor Klee Festschrift (P. Gritzmann and B. Sturmfels, eds.), DIMACS Ser. Discrete Math. Theoret. Comput. Sci., vol. 6, Amer. Math. Soc., Providence, RI, 1991, pp. 237–243.

[17] J. Richter-Gebert, *Realization spaces of 4-polytopes are universal,* Preprint 448/1995, Technical University, Berlin, 1995.

[18] P. Shor, *Stretchability of pseudolines is NP–hard,* Applied Geometry and Discrete Mathematics – The Victor Klee Festschrift (P. Gritzmann, B. Sturmfels, eds.), DIMACS Ser. Discrete Math. Theoret. Comput. Sci., vol. 4, Amer. Math. Soc., Providence, RI, 1991, pp. 531–554.

[19] E. Steinitz, *Polyeder und Raumeinteilungen,* Encyclopädie der Mathematischen Wissenschaften, Band 3 (Geometrie) Teil 3AB12, Teubner, Leipzig, 1922, pp. 1–139.

[20] E. Steinitz and H. Rademacher, *Vorlesungen über die Theorie der Polyeder,* Springer-Verlag, Berlin, 1934; reprint, Springer-Verlag, 1976.

[21] B. Sturmfels, *Boundary complexes of convex polytopes cannot be characterized locally,* Bull. London Math. Soc. **35** (1987), 314–326.

[22] G.M. Ziegler, *Three problems about 4-polytopes,* Polytopes: Abstract, Convex and Computational (T. Bisztriczky, P. McMullen, and A. Weiss, eds.), Kluwer, Dordrecht, 1994, pp. 499–502.

[23] ———, *Lectures on polytopes,* Graduate Texts in Math., vol. 152, Springer-Verlag, New York, 1995.

Technical University Berlin, FB Mathematik, Sekr. 6-1, Strasse des 17. Juni 136, D-10623 Berlin, Germany

*E-mail address*: richter@math.tu-berlin.de

*E-mail address*: ziegler@math.tu-berlin.de